# J. L. DOOB
# (27 FEBRUARY 1910–7 JUNE 2004)

BIOGRAPHICAL NOTE BY M. YOR

*Universités Paris VI et VII and Institut Universitaire de France*

Doob's essential contributions to Probability theory are discussed; this includes the main early results on martingale theory, Doob's *h*-transform, as well as a summary of Doob's three books. Finally, Doob's 'stochastic triangle' is viewed in the light of the stochastic analysis of the eighties.

**1. Biography of J. L. Doob: Some key points.** It is a euphemism to write that the accession of the theory of Probability to the rank of a mathematical discipline took a long time.

After Kolmogorov [25 April 1903 (Tambov)–20 October 1987 (Moscow)], who laid down the axiomatic foundations of Probability Theory, J. L. Doob is, undoubtedly, one of the mathematicians—at first, a specialist in analytic functions—who has worked most for the rigorous mathematization and creation of what was to become, in particular, with the help of his own works, Probability Theory.

Born in 1910 in Cincinnati (Ohio, USA), deceased in 2004 in Urbana (Illinois, USA), J. L. Doob spent his entire career as Professor of Mathematics at the University of Urbana–Champaign (Illinois) between 1935 and 1978. Very much attached to the country in the vicinity of Champaign, he lived there until his death in June 2004.

J. L. Doob wrote, between 1957 and 1963, some very important articles on conditioned Brownian motion, and gave a probabilistic proof to the theorem of Fatou and several extensions of this theorem for the limits at the boundary of (ratios of) harmonic functions. The next paragraphs of this Notice are more particularly dedicated to these works.

J. L. Doob also wrote *three books*, respectively published in 1953, 1984 and 1993, which have had different impacts among probabilists and, more generally, among the mathematical community:

• *The first book* [36]: *Stochastic Processes* (1953), described by P. A. Meyer [M2] as the "Bible of new Probability," exposes, in particular, the theory of continuous-time stochastic processes, for which several measurability and/or regularity properties are put forward in a crucial way. Doob shows how measure theory, which has to be properly developed for this purpose, allows one to solve quite a number of problems in Probability.

One also finds in this book what made it become a success, which is the explanation of martingale theory in discrete and continuous time, as well as—something which some young probabilists may not quite appreciate nowadays!—a detailed study of Itô's stochastic differential equations; indeed, let us emphasize that this book was published in 1953 and that the global "recognition" of Itô's calculus would only really take place from 1969 onward, with H. P. McKean's wonderful little book: *Stochastic Integrals*.

- *The second book* [91]: In *Classical Potential Theory and Its Probabilistic Counterpart* (1984, over 800 pages!), J. L. Doob exposes, with a lot of care and pedagogy, the following:

- in a first part, the theory of Newtonian potential, under its classical form as well as with his personal contributions,
- and in a second part, after having detailed the necessary probability theory, he shows, by referring very precisely to the corresponding points in the first part, how a probabilistic argument allows one to find these results again, or he gives them a more purely probabilistic version; thus, the probabilistic counterpart of a harmonic function is a martingale....

In this second monumental volume, J. L. Doob masterly gathered and exploited all his deep knowledge of these two fields.

One may somehow regret[1] that this treatise, essentially due to its size, seems mostly to have frightened the readers for whom it was meant—potentialists and probabilists, in the first place—and who chose instead smaller treatises on these subjects which had just been written within a few years; let me cite here, for example:

M. Rao: *Brownian Motion and Classical Potential Theory*. Math. Inst., Aarhus Univ. (1977),

R. Durrett: *Brownian Motion and Martingales in Analysis*. Wadsworth, Belmont, Calif. (1984),

S. Port and C. Stone: *Brownian Motion and Classical Potential Theory, Probability and Mathematical Statistics*. Academic Press [Harcourt Brace Jovanovich, Publishers], New York–London (1978).

---

[1]A more positive explanation, of this reading flaw..., may be that at the beginning of the 80's probabilists found great interest in stochastic differential geometry on one hand, and in the stochastic calculus of variations (Malliavin calculus) on the other hand, which is clearly shown by the Proceedings of the Durham Conference (Great Britain) in July 1980, organized by D. Williams. Thus, the ball kept rolling....

The very existence of these books, their titles and their contents show how fundamental the contributions of J. L. Doob's works have been in these fields.

- *The third book* [97]: *Measure Theory* (1994) presents—these are the exact words used by Doob himself in the introduction of this book—"how every training analyst" should approach measure theory, including particularly the probabilistic concepts of independence, conditional independence and martingale, as well as the corresponding theorems (of martingale convergence, for example) and their main applications in the field of analysis.

To summarize, Doob's treatises are rigorous, without any extreme excess in formalization, and go straight to the point, the author himself strongly explaining why such or such notion is important, a way of writing which is not very common among mathematicians....

The reader of this article may be interested to look at other biographical references about J. L. Doob, either of a personal nature ([Sn]), or a scientific one ([B, M1, M2]).

**2. A cursory glance at J. L. Doob's work.** If one may attribute—without any great risk of error—the *discovery* of the notion of martingale to Jean Ville (1910–1988), who exposed in his thesis, *Etude critique de la notion de collectif*, Paris (1939), it is J. L. Doob's works that developed martingale *theory*, by establishing the convergence theorems, and several important uses of martingales. I shall only mention two fundamental results in particular:

- The first is the *stopping theorem*, which expresses that the constant expectation property of a martingale, this expectation being taken at all deterministic times $t$, is still valid when $t$ is replaced by any bounded stopping time $T$. The stopping theorem allows one to make explicit the Laplace transform (or the characteristic function)—and, consequently, the law—of many functionals of Brownian motion and more generally of Markov processes. Consequently, the strategy "find the martingale" has systematically been developed since Doob's works, in order to obtain such results.
- A second fundamental result concerns the bound in $L^p$ ($p > 1$) norm of the supremum of a positive sub-martingale $(X_u, u \leq t)$, by a multiple of the $L^p$ norm of $X_t$. Such basic inequalities allow one to estimate the "size" of a sub-martingale; combined with the inequalities of Burkholder–Davis–Gundy comparing the supremum of a martingale and the square root of its quadratic variation, they allow numerous estimations: for example, they play a key role to show the convergence of Picard's iteration method when solving stochastic differential equations with Lipschitzian coefficients.

If a contemporary probabilist is asked, *what are the main contributions of J. L. Doob regarding probability theory?*, the given answer would very likely consist in the two previous results, to which could be added the discovery and use of the $h$-transform notion—often very simply mentioned as Doob's transform—of a Markov process. This transformation applies to the semigroup of this process, which is multiplied "inside" and divided "outside," by $h$, a harmonic function or, more generally, an excessive function for the process. This transformation often very strongly modifies the nature of the original process: thus, with $h(x) = x$, Doob's transform of real valued Brownian motion $X$, written here as $\widetilde{X}$, has the effect of preventing $\widetilde{X}$ from returning to the origin; in an even better way, this conditioning allows this Brownian motion transform $\widetilde{X}$ to escape toward infinity. This new process $\widetilde{X}$ is nothing else but the Euclidian norm of the 3-dimensional Brownian motion, and the study of the couple $(X, \widetilde{X})$ is at the root of several results concerning the real-valued Brownian motion. This particular case of $h$-relation has aroused numerous extensions, and continues to give rise to many researches, in particular, for multidimensional processes.

*Angular limits, Fatou theorem, Martin boundary.* These themas, which are strongly intertwined, represent an essential part of Doob's contribution to the "boundary limit theorems." This fundamental part of Doob's work will help to illustrate precisely how Doob wrote his second treatise [91] (1984):

- page 641, Doob writes: "in Section 1.XII.19, it was shown that, if $h$ is a strictly positive, harmonic function, in a Greenian domain $D$ of $\mathbb{R}^N$ and if $v$ is a positive superharmonic function on $D$, then $v/h$ has a minimal fine limit in nearly every point (relatively to $M_h$, representative measure of $h$, as the integral of a Martin kernel of $D$). Let us now see a probabilistic formulation equivalent to this Theorem 2.X.8 which asserts that $v/h$ admits a limit along some conditional Brownian paths (translation: Doob refers to the $h$-process)."
- Theorem 2.X.8 appears at page 689, just after Doob has developed—from the beginning of Chapter 2.X, page 668, until page 688—the $h$-processes theory.
- The fundamental particular case, where $D$ is the ball of radius $\delta$ in $\mathbb{R}^N$, is presented in pages 691–693.

The deep discussion of this thema in the treatise of 1984 summarizes Doob's fundamental articles, published, in particular, in the "*Bulletin de la SMF*" [45] (1957) and the "*Annales de l'Institut Fourier*" [53] (1959), in which he gives, respectively, a probabilistic proof and a nonprobabilistic one of the generalized Fatou theorem.

**3. Some extensions of J. L. Doob's work.** To summarize J. L. Doob's main contribution to *Probability Theory*, in the Special Volume of the *Illinois Journal of Mathematics* (2007) dedicated to the memory of Doob and edited by D. Burkholder, Michel Emery introduced the image of Doob's *stochastic triangle* having for apexes:

- $S_1$: a filtered probability space;
- $S_2$: all the stopping times on this space;
- $S_3$: the space of martingales on this space.

One of the aims of Paul-André Meyer, who was among the first probabilists to extend Doob's work, was to classify the stopping times of a filtered probability space.

In order to do so, he had to introduce three fundamental families of processes associated to such a space:

- *predictable processes*, which are "strictly" in the past of the filtration;
- *optional processes*, which represent past and present of the filtration;
- *progressively measurable processes*, only just adapted to the ambient filtration.

Meyer's classification work consisted—among other aims—in describing these stopping times. In particular, all stopping times of a filtered probability space are predictable if and only if all the martingales of this space are continuous,[2] which then gives an illustration of the idea of "Doob's stochastic triangle."

Along with the development of Itô's stochastic calculus, which allows one to integrate every predictable bounded process with respect to a martingale, it was natural to try to unify this calculus with the one of Lebesgue–Stieltjes' integrals, and for this purpose, quite some interest was devoted to the study of *semi-martingales*, that is, the sums of a martingale and a bounded variation process, adapted to the ambient filtration. A famous theorem of Bichteler–Dellacherie characterizes semi-martingales as being the "right integrators" of bounded predictable processes, a very satisfying result as to the nature of the stochastic integral.

Another attempt to extend the stochastic triangle consists in asking oneself which are the random times $\rho$—which are not necessarily stopping times any longer—such that, when a martingale is stopped in $\rho$, one still obtains at least a semi-martingale, in a (new) suitable filtration of course, for which $\rho$ would now be a stopping time.

This question was simultaneously asked—and independently, I think—by D. Williams and P. A. Meyer in 1976–1977; it was partly solved by M. Barlow and T. Jeulin, among others, who showed that if $\rho$ is the end of a predictable set in the original filtration, the question admits a positive answer.

---

[2]This is the case for the Brownian motion filtered probability space.

This was the starting point for the theory of enlargement of filtration, the aim of which is to determine for which super-filtrations of a given filtration every original martingale stays a semi-martingale.

Currently, one may say that two types of enlargements have been developed:

- *initial enlargements* where all new information is brought to the time origin,
- *predictable enlargements* where the new information is brought as time goes by.

The interested reader shall find developments of these points in the following: R. Mansuy, M. Yor: *Random Times and Enlargements of Filtrations in a Brownian Setting*, *LNM* **1873**, Springer (2006).

I just mentioned, in the previous points, how the apexes $S_2$ and $S_3$ of Doob's stochastic triangle were the objects of important developments; this was also the case, even in a previous period actually, for the apex $S_1$, where one tried to understand what becomes of a martingale relative to a filtered probability space, when one modifies, in an absolutely continuous way, the reference probability. There again, original martingales are transformed into semi-martingales for the new filtered probabilistic space; this was described with the help of the various versions of Girsanov's theorem, established under more and more general conditions, originally in a Markovian frame and, finally, simply in the general frame of $S_1$, by Van Schuppen–Wong (1974).

Let us note that in *Financial Mathematics*, the following inverse problem is very important and was solved for a large class of apexes $S_1$: given a *semi-martingale* $X$ relative to a filtered probability space, where $P$ is the reference probability, find all the probabilities $Q$, locally equivalent to $P$, such that under $Q$, $X$ is a martingale. See, for example, the seminal paper of Harrison–Kreps in the *Journal of Economic Theory* (1979), and the fundamental paper of Delbaen–Schachermayer in *Mathematische Annalen* (1994).

Thus, each apex of Doob's stochastic triangle has been the object of several developments, which often necessitated important efforts, and of which in many ways one can say that they are far from being entirely completed, showing thus the depth and originality of J. L. Doob's work. For further discussions of J. L. Doob's work and personality, we refer the reader to [B, M1, M2, Sn].

*Small glossary of the "Stochastic Triangle.*[3]*"*

---

[3]Written originally for nonprobabilists.

- *Probability space*: a triplet made up of the following:
  - a reference space, often noted $\Omega$;
  - a $\sigma$-algebra $\mathcal{F}$ on $\Omega$, that is, a family of subsets of $\Omega$, having some stability properties;
  - a probability $P$, that is, a "$\sigma$-additive" application which associates to every set $A$ of $\mathcal{F}$ a number $P(A)$, with values on [0,1].
- *Filtration*: an increasing family $(\mathcal{F}_t)_{t\geq 0}$ of sub-$\sigma$-algebras of $\mathcal{F}$; $\mathcal{F}_t$ "mathematisizes" the past until time $t$.
- *Stopping time $T$*: an application of $(\Omega, \mathcal{F})$ with values on $[0, \infty]$ such that, for every $t$ $(T \leq t)$ belongs to $\mathcal{F}_t$.

  A typical example of a stopping time is the first time when a phenomenon happens, which can be "observed from the filtration $(\mathcal{F}_t)$."

  On the other hand, a "last time" is typically not a stopping time.
- *Martingale* (relative to a filtered probability space): a family $(X_t)_{t\geq 0}$ of integrable variables, such that $X_t$ is $\mathcal{F}_t$-measurable, and $E[X_t|\mathcal{F}_s] = X_s$ $(s \leq t)$. Typically, the gain until time $t$ in a "fair game."

*Some articles written by J. L. Doob.*

[37] Semi-martingales and subharmonic functions. *Trans. Amer. Math. Soc.* **77** (1954), 86–121.

[45] Conditional Brownian motion and the boundary limits of harmonic functions. *Bull. Soc. Math. France* **85** (1957), 431–458.

[52] A relativized Fatou theorem. *Proc. Nat. Acad. Sci. USA* **45** (1959), 215–222.

[53] A non-probabilistic proof of the relative Fatou theorem. *Ann. Inst. Fourier Grenoble* **9** (1959), 293–300.

[79] What is a martingale? *Amer. Math. Monthly* (1971), 451–463.

## OTHER BIOGRAPHICAL REFERENCES ABOUT J. L. DOOB

LABORATOIRE DE PROBABILITÉS
ET MODÈLES ALÉATOIRES
UNIVERSITÉS PARIS VI ET VII
AND
INSTITUT UNIVERSITAIRE DE FRANCE
4 PLACE JUSSIEU, CASE 188
F-75252 PARIS CEDEX 05
FRANCE